\documentclass{amsart}

\usepackage{fancyhdr}
\usepackage{lastpage}
\usepackage{stmaryrd,yhmath}

\newcommand{\bdis}{\begin{displaymath}}
\newcommand{\edis}{\end{displaymath}}
\newcommand{\be}{\begin{equation}}
\newcommand{\ee}{\end{equation}}
\newcommand{\mbb}{\mathbb}
\newcommand{\mcal}{\mathcal}

\newcommand{\vp}{\varphi}

\newcommand{\vth}{\vartheta}

\newcommand{\mT}{\mathring{T}}

\newcommand{\zf}{\zeta\left(\frac{1}{2}+it\right)}

\newtheorem{theorem}{Theorem}

\theoremstyle{definition}

\theoremstyle{remark}
\newtheorem{remark}[]{Remark}

\newtheorem*{mydef4}{{\bf Corollary}}

\numberwithin{equation}{section}

\begin{document}

\title{Jacob's ladders and some generalizations of certain Ramachandra's inequality}

\author{Jan Moser}

\address{Department of Mathematical Analysis and Numerical Mathematics, Comenius University, Mlynska Dolina M105, 842 48 Bratislava, SLOVAKIA}

\email{jan.mozer@fmph.uniba.sk}


\begin{abstract}
In this paper we obtain some essential generalizations of certain Ramachandra's inequality, i. e. we obtain new lower estimates for the energies of some
complicated signals generated by the Riemann zeta-function on the critical line.
\end{abstract}
\maketitle

\section{Introduction}

\subsection{}

Let
\bdis
\begin{split}
& Z(t)=e^{i\vth(t)}\zf, \\
& \vth(t)=-\frac t2\ln\pi+\text{Im}\ln\Gamma\left(\frac 14+i\frac t2\right)
\end{split}
\edis
be the signal generated by the Riemann zeta-function on the critical line. Next, let us remind the following formulae (see \cite{3}, (9.1), (9.2))
\be \label{1.1}
\begin{split}
& \tilde{Z}^2(t)=\frac{{\rm d}\vp_1(t)}{{\rm d}t},\ \vp_1(t)=\frac 12\vp(t),\ t\geq T_0[\vp], \\
& \tilde{Z}^2(t)=\frac{Z^2(t)}{2\Phi'_\vp[\vp(t)]}=\frac{\left|\zf\right|^2}{\left\{ 1+\mcal{O}\left(\frac{\ln\ln t}{\ln t}\right)\right\}\ln t}, \\
& t\in [T,T+U],\ U\in \left.\left( 0,\frac{T}{\ln T}\right.\right]
\end{split}
\ee
where $\vp(t)$ is the Jacob's ladder, i. e. the exact solution of the nonlinear integral equation (see \cite{2})
\bdis
\begin{split}
& \int_0^{\mu[x(T)]}Z^2(t)e^{-\frac{2}{x(T)}t}{\rm d}t=\int_0^T Z^2(t){\rm d}t, \\
& \mu(y)\geq 7y\ln y,\ \mu(y)\to y=\vp_\mu(T)=\vp(T).
\end{split}
\edis
The function $\vp_1$ is, of course, also the Jacob's ladder.

\subsection{}

Let (see \cite{4}, (2.1))
\bdis
\begin{split}
& \vp_1^0(t)=t,\ \vp_1^1(t)=\vp_1(t),\ \vp_1^2(t)=\vp_1(\vp_1(t)),\dots , \\
& \vp_1^k(t)=\vp_1(\vp_1^{k-1}(t)),\dots
\end{split}
\edis
where $\vp_1^k(t)$ stands for the $k$-th iteration of the Jacob's ladder $\vp_1(t),\ t\in [T,T+U]$. We have obtained the following formula in our paper \cite{4}
\be \label{1.2}
\begin{split}
& \frac 1U\int_T^{T+U}\prod_{k=0}^n\left|\zeta\left(\frac 12+i\vp_1^k(t)\right)\right|^2{\rm d}t\sim \\
& \sim \prod_{k=0}^n\frac{1}{\vp_1^{k}(T+U)-\vp_1^k(T)}\int_{\vp_1^k(T)}^{\vp_1^k(T+U)}\left|\zf\right|^2{\rm d}t, \\
& U\in \left(\left. 0,\frac{T}{\ln^2T}\right.\right],\quad T\to\infty.
\end{split}
\ee
This formula is motivated by the well-known multiplication formula from the theory of probability
\bdis
M\left(\prod_{k=1}^n X_k\right)=\prod_{k=1}^n M(X_k)
\edis
where $X_k$ are the independent random variables and $M$ is the population mean.

\subsection{}

Next, in the paper \cite{5} we have obtained the following formula in the same direction as (\ref{1.2})
\be \label{1.3}
\begin{split}
& \int_T^{T+U}F[\vp_1^{n+1}(t)]\prod_{k=0}^n\left|\zeta\left(\frac 12+i\vp_1^k(t)\right)\right|^2{\rm d}t\sim \\
& \sim \left\{\int_{\vp_1^{n+1}(T)}^{\vp_1^{n+1}(T+U)}\right\}\ln^{n+1}T, \\
& U\in \left(\left. 0,\frac{T}{\ln^2T}\right.\right],\quad T\to\infty
\end{split}
\ee
for every fixed $n\in\mbb{N}$ and for every Lebesgue-integrable function
\bdis
F(t),\ t\in [\vp_1^{n+1}(T),\vp_1^{n+1}(T+U)],\ F(t)\geq 0\- (\leq 0).
\edis

\subsection{}

Let us remind, finally, the inequality of Ramachandra (see \cite{1}, p. 249)
\be \label{1.4}
\begin{split}
& \int_{T-Y}^{T+Y}\left|\zeta^{(r)}\left(\frac 12+it\right)\right|{\rm d}t> AY\left(\ln Y\right)^{r+\frac 14}, \\
& \ln^cT\leq Y\leq T,\ r=0,1,2,\dots
\end{split}
\ee
for every fixed $r$ ($\zeta^{(r)}$ is the $r$-th derivative of the zeta-function). If we use the simple transformation
\bdis
T-Y\to T,\ T+Y\to T+U
\edis
then we obtain the usual form of (\ref{1.4})
\be \label{1.5}
\int_T^{T+U}\left|\zeta^{(r)}\left(\frac 12+it\right)\right|{\rm d}t> BU\left(\ln U\right)^{r+\frac 14}
\ee
with
\be \label{1.6}
3\ln^cT\leq U\leq 2T.
\ee
In this paper we obtain an essential generalization of the Ramachandra's inequality (\ref{1.5}) in the following cases:
\be \label{1.7}
\begin{split}
& U\in \left[ 3\ln^cT,\frac{T}{\ln T}\right], \\
& U\in \left[ T^{\frac 13+2\epsilon},\frac{T}{\ln^2T}\right].
\end{split}
\ee

\section{Results}

\subsection{}

In the first case of (\ref{1.7}) the following theorem holds true.
\begin{theorem}
If
\be \label{2.1}
\vp_1\left\{\left[\mT,\mathring{\widehat{T+U}}\right]\right\}=[T,T+U]
\ee
then
\be \label{2.2}
\begin{split}
& \int_{\mT}^{\mathring{\widehat{T+U}}}\left|\zeta^{(r)}\left(\frac 12+i\vp_1(t)\right)\right|\left|\zf\right|^2{\rm d}t> \\
& > \frac{B}{2}U(\ln U)^{r+\frac 14}\ln T,\ U\in \left[ 3\ln^cT,\frac{T}{\ln T}\right], \\
& r=0,1,2,\dots ,\ T\to\infty
\end{split}
\ee
for every fixed $r$.
\end{theorem}

Next, in the second case of (\ref{1.7}) the following Theorem holds true.

\begin{theorem}
If
\be \label{2.3}
U\in \left[ T^{\frac 13+2\epsilon},\frac{T}{\ln^2T}\right]
\ee
then
\be \label{2.4}
\begin{split}
& \int_T^{T+U}\left|\zeta^{(r)}\left(\frac 12+i\vp_1^{n+1}(t)\right)\right|\prod_{k=0}^n\left|\zeta\left(\frac 12+i\vp_1^k(t)\right)\right|^2{\rm d}t> \\
& > \frac B2U(\ln U)^{r+\frac 14}(\ln T)^{n+1}, \\
& r=0,1,2,\dots
\end{split}
\ee
for every fixed $r,n$; $n\in \mbb{N}$.
\end{theorem}

If we use the Cauchy's inequality in (\ref{2.4}) by the following way
\bdis
\int_T^{T+U}|\zeta^{(r)}|\prod_{k=0}^n|\zeta|^2{\rm d}t\leq \sqrt{U}\left[\int_T^{T+U}|\zeta^{(r)}|^2\prod_{k=0}^n |\zeta|^4\right]^{1/2}, \dots \\
\edis
then we obtain (for example) the following

\begin{mydef4}
\be \label{2.5}
\begin{split}
& \int_T^{T+U}\left|\zeta^{(r)}\left(\frac 12+i\vp_1^{n+1}(t)\right)\right|^{2^m}
\prod_{k=0}^n\left|\zeta\left(\frac 12+i\vp_1^{k}(t)\right)\right|^{2^{m+1}}{\rm d}t> \\
& > \left(\frac B2\right)^{2^m}U(\ln U)^{2^m(r+1/4)}(\ln T)^{2^m(n+1)}
\end{split}
\ee
for every $m\in\mbb{N}$ and every fixed $r,n$.
\end{mydef4}

\begin{remark}
The essential generalizations of the Ramachandra's inequality are expressed by the formulae (\ref{2.2}), (\ref{2.4}) and (\ref{2.5}) for the corresponding segments. Namely, the lower estimates
of the energies of the following complicated signals
\bdis
\begin{split}
& \sqrt{\left|\zeta^{(r)}\left(\frac 12+i\vp_1^{n+1}(t)\right)\right|}\left|\zf\right|,\ t\in \left[\mT,\mathring{\widehat{T+U}}\right], \\
& \sqrt{\left|\zeta^{(r)}\left(\frac 12+i\vp_1^{n+1}(t)\right)\right|}\prod_{k=0}^n\left|\zeta\left(\frac 12+i\vp_1^{k}(t)\right)\right|,\ U\in \left[ T^{\frac 13+2\epsilon},\frac{T}{\ln^2T}\right], \\
& \left|\zeta^{(r)}\left(\frac 12+i\vp_1^{n+1}(t)\right)\right|^{2^{m-1}}\prod_{k=0}^n\left|\zeta\left(\frac 12+i\vp_1^{k}(t)\right)\right|^{2^m},\ U\in \left[ T^{\frac 13+2\epsilon},\frac{T}{\ln^2T}\right]
\end{split}
\edis
are expressed by the formulae (\ref{2.2}), (\ref{2.4}) and (\ref{2.5}).
\end{remark}

\section{Proof of Theorems}

\subsection{}

In the case (\ref{2.1}) we apply our lemma (see \cite{3}, Lemma 4, (9.7)): if
\bdis
f(x)\geq 0 \- (\leq 0),\ x\in [T,T+U]
\edis
is a Lebesgue integrable function then
\be \label{3.1}
\begin{split}
& \int_{\mT}^{\mathring{\widehat{T+U}}} f[\vp_1(t)]\left|\zf\right|^2{\rm d}t= \\
& = \left\{ 1+\mcal{O}\left(\frac{\ln\ln T}{\ln T}\right)\right\}\ln T\int_T^{T+U} f(x){\rm d}x,\ U\in \left(\left. 0,\frac{T}{\ln T}\right.\right].
\end{split}
\ee
If
\bdis
f(x)=\left|\zeta^{(r)}\left(\frac 12+ix\right)\right|
\edis
then we obtain from (\ref{3.1}) by (\ref{1.5})
\bdis
\begin{split}
& \int_{\mT}^{\mathring{\widehat{T+U}}} \left|\zeta^{(r)}\left(\frac 12+i\vp_1(t)\right)\right|\left|\zf\right|^2{\rm d}t= \\
& =\left\{ 1+\mcal{O}\left(\frac{\ln\ln T}{\ln T}\right)\right\}\ln T\int_T^{T+U}\left|\zeta^{(r)}\left(\frac 12+it\right)\right|{\rm d}t> \\
& >\left\{ 1+\mcal{O}\left(\frac{\ln\ln T}{\ln T}\right)\right\} BU(\ln U)^{r+\frac 14}\ln T>\frac B2 U(\ln U)^{r+\frac 14}\ln T,
\end{split}
\edis
and (see (\ref{1.6}), (\ref{3.1}))
\bdis
\left[ 3\ln^cT,2T\right]\bigcap \left(\left. 0,\frac{T}{\ln T}\right.\right]=\left[ 3\ln^cT,\frac{T}{\ln T}\right],
\edis
i. e. we have obtained the estimate (\ref{2.2}).

\subsection{}

In order to prove the inequality (\ref{2.4}) we use the formula (\ref{1.3}) in the case
\bdis
F(t)=\left|\zeta^{(r)}\left(\frac 12+i\vp_1^{n+1}(t)\right)\right|.
\edis
Consequently, we obtain (see (\ref{1.5})) that
\be \label{3.2}
\begin{split}
& \int_T^{T+U}\left|\zeta^{(r)}\left(\frac 12+i\vp_1^{n+1}(t)\right)\right|\prod_{k=0}^n \left|\zeta\left(\frac 12+i\vp_1^{k}(t)\right)\right|{\rm d}t\sim \\
& \sim \left\{\int_{\vp_1^{n+1}(T)}^{\vp_1^{n+1}(T+U)}\left|\zeta^{(r)}\left(\frac 12+it\right)\right|{\rm d}t\right\}(\ln T)^{n+1}> \\
& > (1-\epsilon)B\left[\vp_1^{n+1}(T+U)-\vp_1^{n+1}(T)\right]\left\{\ln\left[\vp_1^{n+1}(T+U)-\vp_1^{n+1}(T)\right]\right\}^{r+\frac 14}(\ln T)^{n+1}.
\end{split}
\ee
In the case (\ref{2.3}), i. e. in the \emph{macroscopic case} with respect to the terminology used in \cite{5}, we have
\be \label{3.3}
\vp_1^{n+1}(T+U)-\vp_1^{n+1}(T)\sim U,\ U\in \left[ T^{\frac 13+2\epsilon},\frac{T}{\ln^2 T}\right],
\ee
(see \cite{5}, (3.8)). Then, finally, from (\ref{3.2}) by (\ref{3.3}) the estimate
\bdis
> \frac B2 U(\ln U)^{r+\frac 14}(\ln T)^{n+1}
\edis
follows, i. e. the formula (\ref{2.4}) is verified.

\section{Concluding formulae}

In our paper \cite{5}, (2.14) -- (2.17), we have obtained the following formulae
\be \label{4.1}
\int_T^{T+U}\prod_{k=0}^n \left|\zeta\left(\frac 12+i\vp_1^{k}(t)\right)\right|{\rm d}t\sim U\ln^{n+1}T,\ U\in \left[ T^{\frac 13+2\epsilon},\frac{T}{\ln^2 T}\right],
\ee
\be \label{4.2}
\begin{split}
& \int_T^{T+U_1} \left|\zeta\left(\frac 12+i\vp_1^{n+1}(t)\right)\right|^4\prod_{k=0}^n \left|\zeta\left(\frac 12+i\vp_1^{k}(t)\right)\right|^2{\rm d}t\sim \\
& \sim \frac{1}{2\pi^2}U_1\ln^{n+5}T,\ U_1=T^{\frac 78+\epsilon},
\end{split}
\ee
\be \label{4.3}
\begin{split}
& \int_T^{T+U}\left\{\arg\zeta\left(\frac 12+i\vp_1^{n+1}(t)\right)\right\}^{2l}\prod_{k=0}^n \left|\zeta\left(\frac 12+i\vp_1^{k}(t)\right)\right|^2{\rm d}t\sim \\
& \sim \frac{(2l)!}{l!4^l}U\ln^{n+1}T(\ln\ln T)^l,\ U\in \left[ T^{\frac 12+\epsilon},\frac{T}{\ln^2 T}\right],
\end{split}
\ee
\be \label{4.4}
\begin{split}
& \int_T^{T+U}\left\{ S_1[\vp_1^{n+1}(t)]\right\}^{2l}\prod_{k=0}^n \left|\zeta\left(\frac 12+i\vp_1^{k}(t)\right)\right|^2{\rm d}t\sim \\
& \sim d_lU\ln^{n+1}T,\ U\in \left[ T^{\frac 12+\epsilon},\frac{T}{\ln^2 T}\right],
\end{split}
\ee
$T\to\infty$, for every fixed $l,n\in\mbb{N}$ where
\bdis
S_1(T)=\int_0^T S(t){\rm d}t,\quad S(t)=\frac{1}{\pi}\arg\zf,
\edis
and the function $\arg$ is defined by the usual way (comp. \cite{6}, p. 179).

\begin{remark}
Now, we obtain from (\ref{4.1}) -- (\ref{4.4}) in the direction of our Corollary the following inequalities
\bdis
\begin{split}
& \int_T^{T+U}\prod_{k=0}^n \left|\zeta\left(\frac 12+i\vp_1^{k}(t)\right)\right|^{2^{m+1}}{\rm d}t>D^{2^m}U(\ln T)^{2^m(n+1)}, \\
& U\in \left[ T^{\frac 13+2\epsilon},\frac{T}{\ln^2 T}\right], \\
& \int_T^{T+U_1} \left|\zeta\left(\frac 12+i\vp_1^{n+1}(t)\right)\right|^{2^{m+2}}\prod_{k=0}^n \left|\zeta\left(\frac 12+i\vp_1^{k}(t)\right)\right|^{2^{m+1}}{\rm d}t > \\
& > D^{2^m}U_1(\ln T)^{2^m(n+5)},\ U_1=T^{\frac 78+2\epsilon},\\
& \int_T^{T+U}\left\{\arg\zeta\left(\frac 12+i\vp_1^{n+1}(t)\right)\right\}^{l2^{m+1}}\prod_{k=0}^n \left|\zeta\left(\frac 12+i\vp_1^{k}(t)\right)\right|^{2^{m+1}}{\rm d}t> \\
& > \{ D(l)\}^{2^m}U(\ln T)^{2^m(n+1)}(\ln\ln T)^{l2^m},\ U\in \left[ T^{\frac 12+\epsilon},\frac{T}{\ln^2 T}\right], \\
& \int_T^{T+U}\left\{ S_1[\vp_1^{n+1}(t)]\right\}^{l2^{m+1}}\prod_{k=0}^n \left|\zeta\left(\frac 12+i\vp_1^{k}(t)\right)\right|^{2^{m+1}}{\rm d}t > \\
& > D^{2^m}U(\ln T)^{2^m(n+1)},\ U\in \left[ T^{\frac 12+\epsilon},\frac{T}{\ln^2 T}\right],
\end{split}
\edis
$T\to\infty$, for the energies of the corresponding complicated signals.
\end{remark}

\thanks{I would like to thank Michal Demetrian for his help with electronic version of this paper.}

\end{document}